\documentclass[11pt]{article}
\usepackage[utf8x]{inputenc}
\usepackage{amssymb}
\usepackage{amsthm}
\usepackage{graphicx}

\newcommand{\Rm}{{\mathbb R}}

\newcommand{\eps}{\varepsilon}

\newtheorem{theorem}{Theorem}[section]
\newtheorem{lemma}[theorem]{Lemma}
\newtheorem{proposition}[theorem]{Proposition}

\begin{document}


\title{I-method for Defocusing, Energy-subcritical Nonlinear Wave Equation}
\author{Ruipeng Shen,\\
        Department of Mathematics,\\
        University of Chicago}
\date{May 2012}
\maketitle
\section{Introduction}
In this paper we will consider the following non-linear wave equation in 3-dimensional space
\begin{equation}
\left\{\begin{array}{l} \partial_t^2 u - \Delta u = F (u), \,\,\, (x,t)\in \Rm^3 \times \Rm;\\
u|_{t=0} = u_0 \in {\dot{H}^s}\cap{\dot{H}^{s_p}} (\Rm^3);\\
\partial_t u |_{t=0} = u_1 \in \dot{H}^{s-1}\cap \dot{H}^{s_p -1}(\Rm^3).\end{array}\right. \label{equ}
\end{equation}
Here the non-linear term $F(u)$ and the coefficients $s_p, s$ are given as below
\begin{eqnarray*}
 F(u) &=& - |u|^{p-1} u.\\
 s_p &=& \frac{3}{2} - \frac{2}{p-1}.
\end{eqnarray*}
\[
 s_p < s < 1.
\]
We will assume $p$ is slightly smaller than $5$, which makes $s_p$ slightly smaller than $1$.
\paragraph{The Energy Space}
If $s=1$, in other words the initial data is in the space $\dot{H}^1 \times L^2$, then the following quantity is
called the energy. The energy is a constant for all time as long as the solution still exists.
\[
 E(u, \partial_t u) = \int_{\Rm^3} (\frac{|\nabla u|^2}{2} + \frac{|\partial_t u|^2}{2} + \frac{|u|^{p+1}}{p+1}) dx.
\]
In this case we are able to obtain global existence and well-posedness of the solution using a basic fixed point argument.
In this paper we are trying to make a weaker assumption, namely, $s$ is greater than $s_p$ but smaller than $1$,
which makes it impossible to use the energy above directly. The I-method described
in many earlier articles (Please see \cite{ImethodcubicNLW, troy1}) can solve this problem for $s$ sufficiently close to $1$.
\paragraph{The Introduction of $I$-operator}
Let us define
\begin{equation} \label{def of I}
 \widehat{I u} (\xi) = \eta(\frac{\xi}{N}) \hat{u} (\xi).
\end{equation}
Here $\eta (\xi)$ is a positive, radial and smooth function defined in $\Rm^3$ such that
\begin{equation}
 \eta (\xi) = \left\{\begin{array}{l} 1, \hbox{ if } |\xi| \leq 1;\cr
                      (\frac{1}{|\xi|})^{1-s}, \hbox{ if } |\xi| > 2.
                     \end{array}\right.
\end{equation}
The number $N \gg 1$ will be determined later. By lemma \ref{Proposition A}, The following quantity is finite and called the energy.
\begin{equation} \label{def of energy}
 E(t) = E(Iu(t),\partial_t Iu(t)) = \int_{\Rm^3} (\frac{|\nabla Iu(x,t)|^2}{2} + \frac{|\partial_t Iu(x,t)|^2}{2} + 
\frac{|Iu(x,t)|^{p+1}}{p+1}) dx.
\end{equation}
Note that $Iu$ is no longer a solution of the original equation (\ref{equ}). Thus the conservation law does not hold any more for this energy. Instead we will introduce an Almost Conservation Law later(See \cite{Almostconsdg} for another example of almost conservation law). The following is our main theorem.

\begin{theorem} \emph{\textbf{(I-method)}} Assume $p \in (11/3, 5)$. There exists $s_0 = s_0 (p) \in (s_p,1)$ such that
if $u$ is a solution of (\ref{equ}) with initial data $(u_0,u_1)$ so that
\[
 \|(u_0,u_1)\|_{\dot{H}^s \times \dot{H}^{s-1}} \leq A;
\]
\[
\|(u_0,u_1)\|_{\dot{H}^{s_p} \times \dot{H}^{s_p-1}} \leq A;
\]
and $s > s_0$, then we have
\[
 \sup_{t \in [0,T]} \|(u(t),\partial_t u(t))\|_{\dot{H}^s \times \dot{H}^{s-1}} \leq C(A,s,p) (1 + T^{\beta(s,p)});
\]
\[
 \sup_{t \in [0,T]} \|(u(t),\partial_t u(t))\|_{\dot{H}^{s_p} \times \dot{H}^{s_p-1}} \leq C(A,s,p) (1 + T^{\beta'(s,p)});
\]
as long as the interval $[0,T]$ is in the maximal lifespan of $u$.
The constant $C(A,s,p)$ above depends on $A,s,p$ only; the exponents $\beta$'s depend on $s,p$ only.
\end{theorem}
\paragraph{Remark} The number $s_0 (p)$ can be given explicitly by
\[
 s_0(p) = \frac{2 + (5-p)s_p}{7-p}.
\]
It is trivial to verify
\[
 s_0(p) \geq \frac{3p-7}{2(p-1)},\; \frac{p-3}{2},\; \frac{3p-5}{2p}.
\]
\paragraph{Comparison with \cite{newroy}} Tristan Roy's recent paper \cite{newroy} studies the same wave equation but makes different assumptions on the initial data. In stead of assuming the initial data is in the space $\dot{H}^{s_p} \times \dot{H}^{s_p-1}$, the author considers localized initial data and obtains similar results using the I-method. More precisely, Roy assumes that the initial data $(u_0,u_1)$ is in the closure of $C_c^\infty (B(0,R)) \times C_c^\infty (B(0,R))$ with respect to the $\dot{H}^s \times \dot{H}^{s-1}$ topology. The difference between \cite{newroy} and my work is
\begin{itemize}
 \item Roy's paper improves the upper bound for $\dot{H}^s \times \dot{H}^{s-1}$ norm of the high frequency part of the solution. It grows more slowly at $T^{\sim (1-s)^2}$ for localized data, thanks to the finite speed of propagation. In contrast, the upper bound grows at $T^{\sim (1-s)}$ in my work if $s$ is close to $1$.
 \item My paper imposes weaker assumptions on the initial data. In fact, any localized data described above is also in the space $\dot{H}^{s_p} \times \dot{H}^{s_p-1} (\Rm^3)$ by the Sobolev embedding.
\end{itemize}
\paragraph{Global Existence} The main theorem actually implies that the solution can never
break down in a finite time. Otherwise the $\dot{H}^s \times \dot{H}^{s-1}$ norm will be bounded in
$[0, T_+)$. But this means the local solution with initial data $(u(T_+ - \eps), \partial_t u (T_+ - \eps))$
would exist at least for some time $T_1$, which would not depend on $\eps$.This is a contradiction when $\eps < T_1$.
We also have the following theorem using a fixed point argument.
\begin{theorem}
Let $s > s_0 (p)$ and assume that $u$ is a solution of (\ref{equ}) with initial data
\[
 (u_0, u_1) \in (\dot{H}^s \cap \dot{H}^{s_p})\times(\dot{H}^{s-1} \cap \dot{H}^{s_p -1}).
\]
If $\|(u_{0,n}-u_0, u_{1,n} - u_1)\|_{\dot{H}^{s_p} \times \dot{H}^{s_p -1}} \rightarrow 0$,
then we have the following limit holds for any given time $t$,
\[
 \|(u_n (t) -u(t) ,\partial_t  u_n(t) - \partial_t u(t))\|_{\dot{H}^{s_p} \times \dot{H}^{s_p -1}} \rightarrow 0.
\]
Here $u_n(t)$ is the solution of (\ref{equ}) with initial data $(u_{0,n}, u_{1,n})$.
\end{theorem}
\section{Preliminary Results}
Local existence and well-posedness of this kind of equations depends on the following Strichartz estimates.

\begin{figure}[h]
\centering
\includegraphics{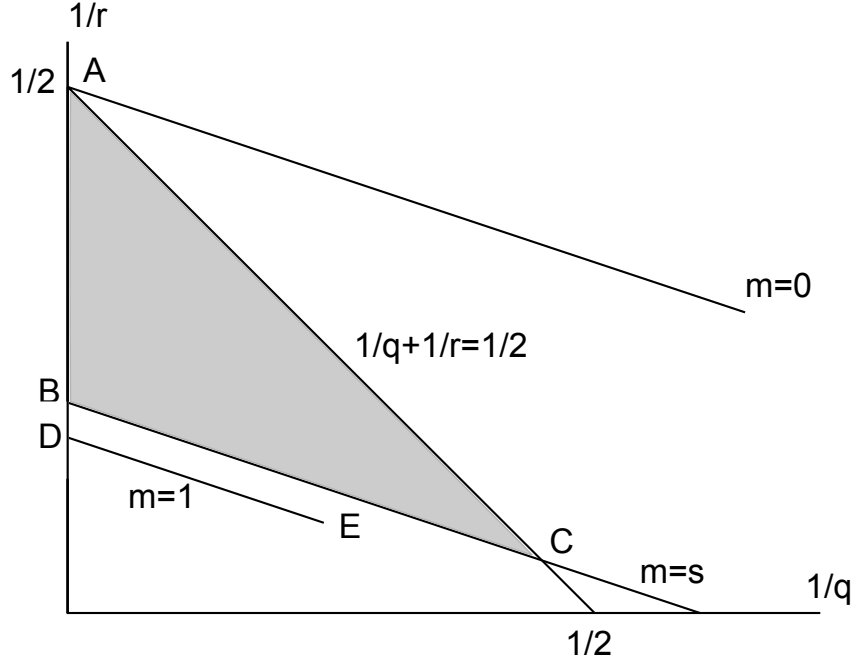}
\caption{Region of allowed pairs}\label{drawing2}
\end{figure}

\begin{proposition} \emph{\textbf{(Generalized Strichartz Inequalities)}}. (Please see proposition 3.1 of \cite{strichartz}, here we use the Sobolev version in $\Rm^3$)
Let $2 \leq q_1,q_2 \leq \infty$, $2 \leq r_1, r_2 < \infty$ and $\rho_1, \rho_2, s \in \Rm$ with
\[
 1/{q_i} + 1/{r_i} \leq 1/2; \,\, i=1,2.
\]
\[
 1/{q_1} + 3/{r_1} = 3/2 - s + \rho_1.
\]
\[
 1/{q_2} + 3/{r_2} = 1/2 + s + \rho_2.
\]
Let $u$ be the solution of the following linear wave equation
\begin{equation}
\left\{\begin{array}{l} \partial_t^2 u - \Delta u = F (x,t), \,\,\,\, (x,t)\in \Rm^3 \times \Rm;\\
u|_{t=0} \in \dot{H}^s (\Rm^3);\\
\partial_t u |_{t=0} = u_1 \in \dot{H}^{s-1}(\Rm^3).\end{array}\right.
\end{equation}
Then we have
\begin{eqnarray*}
 \lefteqn{ \|(u(T), \partial_t u(T))\|_{\displaystyle \dot{H}^s \times \dot{H}^{s-1}} + \|D_x^{\rho_1} u\|_{\displaystyle L^{q_1} L^{r_1} ([0,T]\times \Rm^3)}}\\
 &\leq&
 C \left[ \|(u_0,u_1)\|_{\displaystyle \dot{H}^s \times \dot{H}^{s-1}} + \|D_x^{-\rho_2} F(x,t)\|_{\displaystyle L^{\bar{q}_2} L^{\bar{r}_2}([0,T]\times \Rm^3)}\right].
\end{eqnarray*}
The constant $C$ does not depend on $T$.
\end{proposition}
\paragraph{Remark} In particular, we say that $(q,r)$ is an $m$-admissible pair if
\[
 (\rho_1, s, q_1, r_1) = (0, m, q, r)
\]
satisfies the conditions listed above.
\paragraph{Definition of $Z(J,u)$} Let us assume $p > 11/3$.
In order to take advantage of the Strichartz estimates, we define the following norms
\begin{equation} \label{def1 of z}
 Z_{m,q,r} (J,u) = \|D^{1-m} I u\|_{L_J^q L_x^r}.
\end{equation}
Here $J$ is a closed interval inside the maximum lifespan of the solution $u$, $0 \leq m \leq 1$. The pair $(q,r)$
is $m$-admissible. The Strichartz estimates and the following property of the operator $D^{1-s}I$
\[
 \|D^{1-s} I\|_{L^r \rightarrow L^r} \lesssim 1
\]
show that $Z_{m,q,r} (J,u)$ is always finite if $u$ is a solution of (\ref{equ}). Next step we define
\begin{equation}
 Z(J,u) = \sup_{m,q,r} Z_{m,q,r} (J,u).
\end{equation}
Here the sup is taken among all possible triples $(m,q,r)$ satisfying\\
(I) the pair $(q,r)$ is always $m$-admissible;\\
(II) either
\[
 0 \leq m \leq s,
\]
or
\[
 m=1, \;\hbox{and}\; 1/q \leq \max\left\{ \frac{p-3}{2(p-1)},
\frac{7-p}{4(p-1)} + \frac{1-s}{2(p-1)}\right\}< \frac{1}{2}.
\]
The figure \ref{drawing2} shows all possible pairs $(1/q, 1/r)$ that satisfy the conditions above. This compact region consists of a solid triangle ABC and a closed line segment DE. 
\section{The Proof of Main Theorem}
In this section, we will prove the main theorem. It depends on the following results.

\begin{lemma}\label{Proposition A}
If $(u_0, u_1) \in (\dot{H}^s \cap \dot{H}^{s_p}) \times (\dot{H}^{s-1} \cap \dot{H}^{s_p -1})$, then we have
\[
 \|\nabla I u_0\|_{L^2} \lesssim N^{1-s} \|u_0\|_{\dot{H}^s}.
\]
\[
 \|I u_1\|_{L^2} \lesssim N^{1-s} \|u_1\|_{\dot{H}^{s-1}}.
\]
\[
 \|I u_0\|_{L^{p+1}}^{p+1} \lesssim N^{2(1-s)} \|u_0\|_{\dot{H}^s}^2 \|u_0\|_{\dot{H}^{s_p}}^{p-1}.
\]
In summary, we have
\begin{equation} \label{es of e}
 E(I u_0, I u_1) \lesssim N^{2(1-s)} \left( \|u_0\|_{\dot{H}^s}^2 + \|u_1\|_{\dot{H}^{s-1}}^2 +
 \|u_0\|_{\dot{H}^s}^2 \|u_0\|_{\dot{H}^{s_p}}^{p-1} \right).
\end{equation}
\end{lemma}
\begin{lemma}\label{Proposition B}
Let $u$ be a solution of the equation (\ref{equ}), $J=[0,T]$, then
\begin{eqnarray*}
 \lefteqn{\|(u(T), \partial_t u(T))\|_{{\dot{H}^s}\times{\dot{H}^{s-1}}} \leq
 \| (u(0), \partial_t u(0))\|_{{\dot{H}^s}\times{\dot{H}^{s-1}}}\label{prop c}}\\
 && + C_{s,p} \left(\sup_{t \in J} E(t)^{1/2} +  T \sup_{t \in J} E(t)^{\frac{p}{p+1}} +
\frac{Z^p(J,u)}{N^{\frac{5-p}{2} + 1-s}}\right).
\end{eqnarray*}
\end{lemma}

\begin{lemma}\label{Proposition C}
There exist $\tau_0 = \tau_0 (s, p)$ and $N_0 = N_0 (s, p)$ such that if
$|J| \leq \tau_0$, $N > N_0$ and $u(x,t)$ is a solution of the equation with
\[
 \sup_{t \in J} E(t) \leq 1,
\]
then we have
\[
 Z(J,u) \lesssim_{s,p} 1.
\]
\end{lemma}
\begin{lemma}\emph{\textbf{Almost Conservation Law of Energy}}\label{almost conservation law}
If $u$ is a solution of the equation, then the inequality
\[
 |E(Iu(t_1)) - E(Iu(t_2))| \lesssim \sup_{t \in J} E(t)^{1/2} \frac{Z^p (J,u)}{N^{(5-p)/2}}
\]
holds for all times $t_1, t_2 \in J$.
\end{lemma}
\noindent These lemmas will be proved in the later sections. Now let us show that the main
theorem holds assuming these lemmas.

\paragraph{Step 1: Scaling} Let $u_\lambda$ be
\begin{equation} \label{def of scaling}
 u_{\lambda}(x,t) = \frac{1}{\lambda^{\frac{3}{2}- s_p}} u (\frac{x}{\lambda} , \frac{t}{\lambda}).
\end{equation}
Thus
\begin{equation}
 \partial_t u_{\lambda}(x,t) = \frac{1}{\lambda^{\frac{5}{2}- s_p}} \partial_t u (\frac{x}{\lambda} , \frac{t}{\lambda}).
\end{equation}
If $u(x,t)$ is a solution of the equation (\ref{equ}), one can check that $u_\lambda$ is still a solution of the original equation. In addition,
the $\dot{H}^{s_p} \times \dot{H}^{s_p -1}$ norm is preserved under this rescaling.
Using (\ref{es of e}) we know the energy
\begin{eqnarray*}
 && E(I u_\lambda(0), I \partial_t u_\lambda (0))\\
 &\lesssim& N^{2(1-s)} \left( \|u_\lambda (0)\|_{\dot{H}^s}^2 + \|\partial_t u_\lambda (0)\|_{\dot{H}^{s-1}}^2 +
 \|u_\lambda (0)\|_{\dot{H}^s}^2 \|u_\lambda (0)\|_{\dot{H}^{s_p}}^{p-1} \right)\\
 &\lesssim&  N^{2(1-s)} \lambda^{ 2s_p - 2s} \left( \|u (0)\|_{\dot{H}^s}^2 + \|\partial_t u (0)\|_{\dot{H}^{s-1}}^2 +
 \|u (0)\|_{\dot{H}^s}^2 \|u (0)\|_{\dot{H}^{s_p}}^{p-1} \right)\\
 &=&  N^{2(1-s)} \lambda^{ 2s_p - 2s} \left( \|u_0\|_{\dot{H}^s}^2 + \|u_1\|_{\dot{H}^{s-1}}^2 +
 \|u_0\|_{\dot{H}^s}^2 \|u_0\|_{\dot{H}^{s_p}}^{p-1} \right)
\end{eqnarray*}
Let us define $C(u)$ by
\begin{equation}
 C(u) =  \|u_0\|_{\dot{H}^s}^2 + \|u_1\|_{\dot{H}^{s-1}}^2 +
 \|u_0\|_{\dot{H}^s}^2 \|u_0\|_{\dot{H}^{s_p}}^{p-1},
\end{equation}
and choose
\begin{equation} \label{choice of lambda}
 \lambda = C_{s,p} C(u)^{\frac{1}{2(s -s_p)}} N^{\frac{1 -s}{s- s_p}}.
\end{equation}
If $C_{s,p}$ is sufficiently large, then
\[
  E(I u_\lambda(0), I \partial_t u_\lambda (0)) \leq 1/2.
\]

\paragraph{Step 2} We will show that the energy $ E( I u_\lambda (t), I \partial_t u_\lambda (t)) $ is always less than $3/4$ in the whole interval $[0,\lambda T]$ if we choose sufficiently large $N = N(s,p,\|u\|, T)$. Let us define
\[
 T' = \max \left\{t : t \in [0, \lambda T]\,\hbox{such that}\, E( I u_\lambda (t'), I \partial_t u_\lambda (t')) \leq 3/4,
 \hbox{ for all } t' \in [0,t] \right\}.
\]
By continuity of the energy, if $T' < \lambda T$, we have there exists $\eps >0$, such that
\[
 E( I u_\lambda (t), I \partial_t u_\lambda (t)) \leq 1,
 \hbox{ for all } t \in [0,T' + \eps].
\]
Break the interval $[0, T'+\eps]$ into subintervals $\{J_i\}$, $i=1,2,\cdots, n$, such that
$|J_i| \leq \tau_0$. The constant $\tau_0 = \tau_0 (s,p)$ here and $N_0 (s,p)$ mentioned below are the same constants as in lemma \ref{Proposition C}. We can always choose
\begin{equation}\label{num of subi}
 n \leq \frac{\lambda T}{\tau_0} + 1.
\end{equation}
By lemma \ref{Proposition C}, we have (Let $N > N_0 (s,p)$)
\begin{equation}\label{boundedness of Z}
 Z(J_i, u_\lambda) \lesssim 1.
\end{equation}
Applying Almost Conservation Law in each subinterval, we obtain
\begin{equation}
 E(I u_\lambda(t), I \partial_t u_\lambda (t)) \leq 1/2 + \frac{C_{s,p} \cdot k}{N^{(5-p)/2}}.
\end{equation}
for any $t \in J_k$. Using (\ref{num of subi}) we have for any $t \in [0, T' +\eps]$,
\begin{eqnarray*}
 E(I u_\lambda(t), I \partial_t u_\lambda (t)) &\leq& 1/2 +
\frac{ C_{s,p} (1 + \frac{\lambda T}{\tau_0}) }{N^{(5-p)/2}}\\
 &\leq& 1/2 +
\frac{ C_{s,p} (1 + \lambda T) }{N^{(5-p)/2}}\\
 &\leq& 1/2 +
\frac{ C_{s,p} (1 + C(u)^{\frac{1}{2(s -s_p)}} N^{\frac{1 -s}{s- s_p}} T) }{N^{(5-p)/2}}\\
 &\leq& 1/2 + \frac{C_{s,p}}{N^{\frac{5-p}{2}}} +
\frac{C_{s,p} C(u)^{\frac{1}{2(s -s_p)}} T }{N^{\frac{5-p}{2} - \frac{1-s}{s-s_p}}}.
\end{eqnarray*}
Here we use the choice of $\lambda$ (\ref{choice of lambda}). The constants above $C_{s,p}$ may
be different in each step, but they only depend on the numbers $s,p$. Our assumption on $s$ actually implies
\[
 \frac{5-p}{2} >  \frac{1-s}{s-s_p}.
\]
Choosing
\begin{equation}\label{choice of N}
 N = C \max \left\{C(u)^{\displaystyle \frac{1}{(5-p)(s-s_p) -2(1-s)}}
T^{\displaystyle \frac{1}{\frac{5-p}{2} - \frac{1-s}{s-s_p}}},\, N_0(s,p)\right\},
\end{equation}
we have
\[
 E(I u_\lambda(t), I \partial_t u_\lambda (t)) \leq 3/4
\]
for all $t \in [0,T' + \eps]$ if $C= C(s,p)$ is sufficiently large. This is a contradiction. Thus if we choose $N$ as (\ref{choice of N}),
then the following inequality
\begin{equation}\label{boundedness of E}
 E(I u_\lambda(t), I \partial_t u_\lambda (t)) \leq 3/4
\end{equation}
holds for each $t \in [0,\lambda T]$. Breaking this interval into subintervals $J_i$ as above, we still have
(\ref{boundedness of Z}) and (\ref{num of subi}) holds.
\paragraph{Step 3} Applying lemma \ref{Proposition B} to each subinterval and conducting an induction, we obtain for each $t_0 \in [0,\lambda T]$,
(Use (\ref{num of subi}), (\ref{boundedness of Z}), (\ref{choice of N}) and (\ref{boundedness of E}))
\begin{eqnarray*}
 \lefteqn{\|(u_\lambda (t_0), \partial_t u_\lambda (t_0))\|_{{\dot{H}^s}\times{\dot{H}^{s-1}}}}\\ &\leq&
 \| (u_\lambda (0), \partial_t u_\lambda (0))\|_{{\dot{H}^s}\times{\dot{H}^{s-1}}}
 + C_{s,p} \left( n + \lambda T + \frac{n}{N^{\frac{5-p}{2} + 1-s}}\right)\\
&\leq&
 \lambda^{s_p -s}\| (u (0), \partial_t u (0))\|_{{\dot{H}^s}\times{\dot{H}^{s-1}}}
 + C_{s,p} \left( \lambda T + 1 + \frac{\lambda T  +1}{N^{\frac{5-p}{2} + 1-s}}\right)\\
&\leq&
 \lambda^{s_p -s}\| (u (0), \partial_t u (0))\|_{{\dot{H}^s}\times{\dot{H}^{s-1}}}
 + C_{s,p} \left( \lambda T + 1\right).
\end{eqnarray*}
Rescaling back we have
\begin{eqnarray*}
 \left\|(u (\frac{t_0}{\lambda}), \partial_t u (\frac{t_0}{\lambda}))\right\|_{{\dot{H}^s}\times{\dot{H}^{s-1}}} &\leq&
 \|(u (0), \partial_t u (0))\|_{{\dot{H}^s}\times{\dot{H}^{s-1}}} + C_{s,p} \lambda^{s - s_p} (\lambda T + 1)\\
 &\leq& \|(u_0, u_1)\|_{{\dot{H}^s}\times{\dot{H}^{s-1}}}  + C_{s,p} \lambda^{s - s_p} N^{(5-p)/2}\\
 &\leq& \|(u_0, u_1)\|_{{\dot{H}^s}\times{\dot{H}^{s-1}}} + C_{s,p} C(u)^{1/2} N^{1-s} N^{(5-p)/2}\\
 &\lesssim& \|(u_0, u_1)\|_{{\dot{H}^s}\times{\dot{H}^{s-1}}} + C(u)^{1/2}  N^{1-s + \frac{5-p}{2}}\\
 &\lesssim& \|(u_0, u_1)\|_{{\dot{H}^s}\times{\dot{H}^{s-1}}} + C(u)^\alpha T^{\beta} + C(u)^{1/2}.
\end{eqnarray*}
The exponents $\alpha$ and $\beta$ are given by
\begin{eqnarray}
 \alpha = \alpha(s,p) &=& {\frac{ \frac{5-p}{2}(1 + s -s_p)}{(5-p)(s-s_p) -2(1-s)}};\\
 \beta = \beta(s,p) &=& \frac{1-s + \frac{5-p}{2}}{\frac{5-p}{2} - \frac{1-s}{s-s_p}}.
\end{eqnarray}
In summary
\begin{equation} \label{estimate of hs}
 \sup_{t \in [0,T]} \|(u (t), \partial_t u (t))\|_{{\dot{H}^s}\times{\dot{H}^{s-1}}}
\lesssim \|(u_0, u_1)\|_{{\dot{H}^s}\times{\dot{H}^{s-1}}} +C(u)^{1/2} + C(u)^\alpha T^{\beta}.
\end{equation}
This gives the bound for the $\dot{H}^s \times \dot{H}^{s-1}$ norm. We can also find an upper bound
for the $\dot{H}^{s_p} \times \dot{H}^{s_p -1}$ norm as below.
\paragraph{The $\dot{H}^{s_p} \times \dot{H}^{s_p -1}$ Norm}
By the local theory of the equation with initial data $(u_0,u_1) \in \dot{H}^s \times \dot{H}^{s-1}$,
we know local solution will exist at least in the interval $[0,T_1]$, where the number $T_1$ is given by
\[
 T_1 =  \frac{C_{s,p}}{\|(u_0,u_1)\|_{\dot{H}^s \times \dot{H}^{s-1}}^{\frac{1}{s-s_p}}}.
\]
This is different from the local theory with initial data in the critical space $\dot{H}^{s_p} \times \dot{H}^{s_p -1}$.
In addition, Given each $s$-admissible pair $(q,r)$, we have
\begin{equation}
 \|u\|_{L_t^q L_x^r ([0,T_1] \times \Rm^3)} \lesssim \|(u_0,u_1)\|_{\dot{H}^s \times \dot{H}^{s-1}}.
\end{equation}
Now let the letter $M$ represent the upper bound as below. Please note that we can estimate $M$ by (\ref{estimate of hs}).
\[
 M = \sup_{t \in [0,T]} \|(u (t), \partial_t u (t))\|_{{\dot{H}^s}\times{\dot{H}^{s-1}}}.
\]
If we break the interval $[0,T]$ into subintervals $J_i (i=1,2,\cdots, n)$, such that
\[
|J_i| \leq T_1\approx 1 / M^{\frac{1}{s -s_p}},
\]
then the local theory can be applied in each subinterval. Choosing a specific
$s$-admissible pair $(p/(1 - p(s-s_p)), 6p/(5-2s_p))$, we have
\[
 \|F(u)\|_{L_{J_i}^\frac{1}{1-p(s-s_p)} L_x^{\frac{6}{5-2s_p}}}
\leq \|u\|_{L_{J_i}^\frac{p}{1-p(s-s_p)} L_x^{\frac{6p}{5-2s_p}}}^p \lesssim M^p.
\]
Thus
\[
 \|F(u)\|_{L_{J_i}^1 L_x^{\frac{6}{5-2s_p}}} \leq |J_i|^{p(s-s_p)}
\|F(u)\|_{L_{J_i}^\frac{1}{1-p(s-s_p)} L_x^{\frac{6}{5-2s_p}}}
\lesssim 1.
\]
The bound in question is given by a straightforward computation using the Strichartz estimates as below
\begin{eqnarray*}
 \lefteqn{\|(u,\partial_t u)\|_{C([0,T];\dot{H}^{s_p} \times \dot{H}^{s_p-1})}}\\ &\leq&
 \|(u_0,u_1)\|_{\dot{H}^{s_p} \times \dot{H}^{s_p-1}} + C_{s,p}
\|F(u)\|_{L_{[0,T]}^1 L_x^{\frac{6}{5-2s_p}}} \\
 &\leq&  \|(u_0,u_1)\|_{\dot{H}^{s_p} \times \dot{H}^{s_p-1}} + C_{s,p} (1 + \frac{T}{T_1})\\
 &\leq& \|(u_0,u_1)\|_{\dot{H}^{s_p} \times \dot{H}^{s_p-1}} + C_{s,p} (1 + T M^{\frac{1}{s-s_p}})\\
 &\leq& \|(u_0,u_1)\|_{\dot{H}^{s_p} \times \dot{H}^{s_p-1}} \\
 && + C_{s,p} \left( 1 + T \left(\|(u_0, u_1)\|_{{\dot{H}^s}\times{\dot{H}^{s-1}}} + C(u)^{1/2}
+ C(u)^\alpha T^{\beta}\right)^{\frac{1}{s-s_p}}\right)\\
 &\lesssim& \|(u_0,u_1)\|_{\dot{H}^{s_p} \times \dot{H}^{s_p-1}} + 1 +
T \|(u_0, u_1)\|_{{\dot{H}^s}\times{\dot{H}^{s-1}}}^{\frac{1}{s-s_p}} + T C(u)^{\frac{1}{2(s-s_p)}}\\
&& +
C(u)^{\frac{\alpha}{s-s_p}} T^{\frac{\beta}{s-s_p} + 1}.
\end{eqnarray*}

\section{Proof of Lemma \ref{Proposition A}}
This lemma comes from some basic computation.
\begin{eqnarray*}
 \|\nabla I u_0\|_{L^2}^2 &\lesssim& \int_{|\xi| \leq N} |\xi|^2 |\hat{u_0} (\xi) |^2 d \xi +
\int_{|\xi| > N} |\xi|^2 \frac{N^{2(1-s)}}{|\xi|^{2(1-s)}} |\hat{u_0} (\xi)|^2 d \xi\\
 &\lesssim& \int_{|\xi| \leq N} N^{2(1-s)} |\xi|^{2s} |\hat{u_0} (\xi) |^2 d \xi +
\int_{|\xi| > N}  {N^{2(1-s)}} |\xi|^{2s} |\hat{u_0} (\xi)|^2 d \xi\\
 &\lesssim& N^{2(1-s)} \|u_0\|_{\dot{H}^s}^2.
\end{eqnarray*}
Similar argument shows
\begin{eqnarray*}
 \lefteqn{\|I u_1\|_{L^2}^2 \lesssim \int_{|\xi| \leq N} |\hat{u_1} (\xi) |^2 d \xi +
\int_{|\xi| > N} \frac{N^{2(1-s)}}{|\xi|^{2(1-s)}} |\hat{u_1} (\xi)|^2 d \xi}\\
 &\lesssim& \int_{|\xi| \leq N} N^{2(1-s)} |\xi|^{2(s-1)} |\hat{u_0} (\xi) |^2 d \xi +
\int_{|\xi| > N}  {N^{2(1-s)}} |\xi|^{2(s-1)} |\hat{u_0} (\xi)|^2 d \xi\\
 &\lesssim& N^{2(1-s)} \|u_1\|_{\dot{H}^{s-1}}^2.
\end{eqnarray*}
For the third inequality we have
\begin{eqnarray*}
 \|I u_0\|_{L^{p+1}}^{p+1} &\lesssim& \|I u_0\|_{L^6}^2 \| I u_0 \|_{L^{\frac{3(p-1)}{2}}}^{p-1}\\
  &\lesssim& \|\nabla I u_0\|_{L^2}^2 \| u_0\|_{L^{\frac{3(p-1)}{2}}}^{p-1}\\
  &\lesssim& N^{2(1-s)} \|u_0\|_{\dot{H}^s}^2 \|u_0\|_{\dot{H}^{s_p}}^{p-1}.
\end{eqnarray*}

\section{Proof of Lemma \ref{Proposition B}}
In this section we give the proof of lemma \ref{Proposition B}. We will first estimate the low frequency part, which is more difficult.
By Strichartz estimate we have
\begin{eqnarray*}
 \|P_{\leq 1} (u(T), \partial_t u(T))\|_{{\dot{H}^s}\times{\dot{H}^{s-1}}} &\leq&
 \|S(t) P_{\leq 1} (u(0), \partial_t u(0))\|_{{\dot{H}^s}\times{\dot{H}^{s-1}}}\\  && +
 \left\|\int_0^T \frac{\sin{((T-t)\sqrt{-\Delta})}}{\sqrt{-\Delta}} P_{\leq 1} F(u(t)) dt\right\|_{{\dot{H}^s}\times{\dot{H}^{s-1}}}\\
 &\leq& \|P_{\leq 1} (u(0), \partial_t u(0))\|_{{\dot{H}^s}\times{\dot{H}^{s-1}}} +
 C_{s} \|P_{\leq 1} F(u)\|_{L_J^1 L_x^{\frac{6}{5-2s}}}.
\end{eqnarray*}
We can break the nonlinear part into
\begin{eqnarray*}
 \|P_{\leq 1} F(u)\|_{L_J^1 L_x^{\frac{6}{5-2s}}} &\leq&
\|P_{\leq 1} F(I u)\|_{L_J^1 L_x^{\frac{6}{5-2s}}} +
\|P_{\leq 1} [F(u) - F(I u)]\|_{L_J^1 L_x^{\frac{6}{5-2s}}}\\
 &\lesssim& \|P_{\leq 1} F(I u)\|_{L_J^1 L_x^{\frac{6}{5-2s}}} +
\|F(u) - F(I u)\|_{L_J^1 L_x^{\frac{6}{5-2s}}}\\
 &=& X_1 + X_2,
\end{eqnarray*}
and deal with each part individually
\begin{eqnarray*}
 X_1 &\lesssim& \|P_{\leq 1} F(I u)\|_{L_J^1 L_x^{\frac{p+1}{p}}}\\
 &\lesssim& \|F(I u)\|_{L_J^1 L_x^{\frac{p+1}{p}}}\\
 &\lesssim& T \|F(Iu)\|_{L_J^\infty L_x^{\frac{p+1}{p}}}\\
 &\lesssim& T \sup_{t \in J} \left(\int_{\Rm^3} |F(Iu)|^{\frac{p+1}{p}} dx\right)^{\frac{p}{p+1}}\\
 &\lesssim& T \sup_{t \in J} \left(\int_{\Rm^3} |Iu|^{p+1} dx\right)^{\frac{p}{p+1}}\\
 &\lesssim& T \sup_{t \in J} E(t)^{\frac{p}{p+1}},
\end{eqnarray*}
and (Similar argument is used in the proof of almost conservation law)
\begin{eqnarray*}
 X_2 &\lesssim& \|(|Iu| + |u|)^{p-1} |Iu -u|\|_{L_J^1 L_x^{\frac{6}{5-2s}}}\\
 &\lesssim& \|P_{\ll N} u\|_{L_J^{\frac{1}{\frac{7-p}{4(p-1)} + \frac{1-s}{2(p-1)}}}
 L_x^{\frac{1}{\frac{p-3}{4(p-1)} - \frac{1-s}{6(p-1)}}}}^{p-1} \cdot
\| P_{\gtrsim N} u \|_{L_J^\frac{1}{\frac{p-3}{4} - \frac{1-s}{2}}
L_x^{\frac{1}{\frac{5-p}{4}+ \frac{1-s}{2}}}}\\
&& + \| P_{\gtrsim N} u \|_{L_J^{p} L_x^{\frac{6p}{5-2s}}}^{p-1} \cdot
\|P_{\gtrsim N} u\|_{L_J^p L_x^{\frac{6p}{5-2s}}}\\
 &\lesssim& \|D^{1-1} I u\|_{L_J^{\frac{1}{\frac{7-p}{4(p-1)} + \frac{1-s}{2(p-1)}}}
 L_x^{\frac{1}{\frac{p-3}{4(p-1)} - \frac{1-s}{6(p-1)}}}}^{p-1} \\
&& \times \frac{1}{N^{\frac{5-p}{2} + 1-s}} \| D^{1 - \left(\frac{p-3}{2}-(1-s)\right)} I u \|_{L_J^\frac{1}{\frac{p-3}{4} - \frac{1-s}{2}}
L_x^{\frac{1}{\frac{5-p}{4}+ \frac{1-s}{2}}}}\\
&& + \frac{1}{N^{\frac{5-p}{2} + 1-s}} \| D^{1 - \frac{3p -7 +2s}{2p}} I u \|_{L_J^{p} L_x^{\frac{6p}{5-2s}}}^{p}\\
&\lesssim& \frac{Z^p(J,u)}{N^{\frac{5-p}{2} + 1-s}}.
\end{eqnarray*}
Thus in summary we have
\begin{eqnarray*}
 \|P_{\leq 1} (u(t), \partial_t u(t))\|_{{\dot{H}^s}\times{\dot{H}^{s-1}}} &\leq&
 \|P_{\leq 1} (u(0), \partial_t u(0))\|_{{\dot{H}^s}\times{\dot{H}^{s-1}}}\\
 && + C_{s,p} \left( T \sup_{t \in J} E(t)^{\frac{p}{p+1}} +
\frac{Z^p(J,u)}{N^{\frac{5-p}{2} + 1-s}}\right).
\end{eqnarray*}
Next let us consider the high frequency part
\begin{eqnarray*}
 \|P_{>1} u(t)\|_{\dot{H}^s}^2 &\leq&
 \int_{1 < |\xi| \leq 2N} |\xi|^{2s} |\hat{u}(t,\xi)|^2 d \xi +
 \int_{|\xi| > 2N} |\xi|^{2s} |\hat{u}(t,\xi)|^2 d \xi \\
 &\lesssim& X_1 + X_2.
\end{eqnarray*}
These two terms can be dominated by the energy just at the time $t$.
\[
 X_1 \lesssim \int_{1 < |\xi| \leq 2N} |\xi|^2 |\hat{u}(t,\xi)|^2 d \xi
 \lesssim \|\nabla Iu(t)\|_{L^2}^2 \lesssim E(t).
\]
\begin{eqnarray*}
 X_2 &\lesssim& \frac{1}{N^{2(1-s)}}
\int_{|\xi| > 2N} |\xi|^2 \frac{N^{2(1-s)}}{|\xi|^{2(1-s)}} |\hat{u}(t,\xi)|^2 d \xi\\
 &\lesssim& \frac{1}{N^{2(1-s)}} \|\nabla I u(t)\|_{L^2}^2\\
 &\lesssim& \frac{1}{N^{2(1-s)}} E(t).
\end{eqnarray*}
By similar argument we can show
\[
 \|P_{>1} \partial_t u(t)\|_{\dot{H}^{s-1}}^2 \lesssim E(t).
\]
Combining the low and high frequency parts, we have
\begin{eqnarray*}
 \lefteqn{\|(u(T), \partial_t u(T))\|_{{\dot{H}^s}\times{\dot{H}^{s-1}}} \leq
 \| (u(0), \partial_t u(0))\|_{{\dot{H}^s}\times{\dot{H}^{s-1}}}}\\
 && + C_{s,p} \left(\sup_{t \in J} E(t)^{1/2} +  T \sup_{t \in J} E(t)^{\frac{p}{p+1}} +
\frac{Z(J,u)^p}{N^{\frac{5-p}{2} + 1-s}}\right).
\end{eqnarray*}
\section{Proof of Lemma \ref{Proposition C}}
In this section we will prove lemma \ref{Proposition C}.
\paragraph{Step 1} Let us first consider the estimate for $q = \infty$. Using the Sobolev embedding, we have
\[
 \|D^{1-m} I u\|_{L_J^\infty L_x^r} \lesssim \|D I u\|_{L_J^\infty L_x^2} \lesssim \sup_{t \in J} (E(t))^{1/2} \leq 1.
\]
Thus the estimate holds for $q = \infty$.
\paragraph{Step 2} Now we will first establish an estimate for $m \leq s$. WLOG, let $J = [0,\tau]$.
Applying the operator $D^{1-m} I$ to the original equation (\ref{equ}) and
then using the Strichartz estimate, we obtain
\begin{eqnarray*}
 \|D^{1-m} I u\|_{L_J^q L_x^r} &\lesssim& \|(D^{1-m} I u(0),D^{1-m} I \partial_t u(0))\|_{\dot{H}^m \times \dot{H}^{m-1}}\nonumber \\
 && + \|D^{1-m} I F(u)\|_{L_J^1 L_x^{\frac{6}{5-2m}}}.\label{strichartz01}
\end{eqnarray*}
Using the fact $m \leq s$ we have
\begin{eqnarray*}
 Z_{m,q,r} &\lesssim& \|(I u(0),\partial_t I u(0))\|_{\dot{H}^{1} \times L^2} +
\|D^{1-m} I u\|_{L_J^\infty L_x^\frac{6}{3-2m}} \|u\|_{L_J^{p-1} L_x^{3(p-1)}}^{p-1}\\
 &\lesssim& \sup_{t \in J} E(t)^{1/2}\\
&+& Z_{m, \infty, \frac{6}{3-2m}}
 \left( \tau^{\frac{5-p}{2}} \|P_{\ll N} u\|_{L_J^{\frac{2(p-1)}{p-3}} L_x^{3(p-1)}}^{p-1}
+ \|P_{\gtrsim N} u\|_{L_J^{p-1} L_x^{3(p-1)}}^{p-1} \right)\\
 &\lesssim& \sup_{t\in J} E(t)^{1/2} +
\tau^{\frac{5-p}{2}} \|D^{1-1} I u\|_{L_J^{\frac{2(p-1)}{p-3}} L_x^{3(p-1)}}^{p-1}\\
&+& \frac{1}{N^{\frac{5-p}{2}}}\|D^{1-{s_p}} I u\|_{L_J^{p-1} L_x^{3(p-1)}}^{p-1}.
\end{eqnarray*}
We also need to estimate the case when $m=1$. In this case we have
\begin{eqnarray*}
 \|Iu\|_{L_J^q L_x^r} &\lesssim& \|(I u(0), \partial_t I u(0))\|_{\dot{H}^1 \times L^2} +
\|I F(u)\|_{L_J^1 L_x^2}\\
 &\lesssim& \sup_{t \in J} E(t)^{1/2} + \|D^{1-s} I F(u)\|_{L_J^1 L_x^{\frac{6}{5-2s}}}.
\end{eqnarray*}
Using the same argument as the case $m=s$, we can find the same upper bound as the previous case. In summary
\begin{equation}
 Z(J,u) \leq C_{s,p} \left(\begin{array}{r} \sup_{t\in J} E(t)^{1/2}
 + \tau^{\frac{5-p}{2}} \|D^{1-1} I u\|_{L_J^{\frac{2(p-1)}{p-3}} L_x^{3(p-1)}}^{p-1}\\
+ \frac{1}{N^{\frac{5-p}{2}}}\|D^{1-s_p} I u\|_{L_J^{p-1} L_x^{3(p-1)}}^{p-1}\end{array}\right)\label{st2}
\end{equation}
\paragraph{Remark} It seems that the constant $C_{s,p}$ should have depended on $q,r$ besides $s,p$, because the best constant in a Srtrichartz
estimate depends on the coefficients $(q,r)$. However, it is still possible to find a universal constant that works for each allowed triple.
We can first establish individual estimates as above for those $(1/q,1/r)$ that respond to the vertices (A,B,C,D,E) in the figure \ref{drawing2} and then use an interpolation to gain a universal constant $C_{s,p}$ for all possible triples.
\paragraph{Step 3} Let
\[
 \tilde{Z}(t,u) = \max\left\{ \|D^{1-1} I u\|_{L_{[0,t]}^{\frac{2(p-1)}{p-3}} L_x^{3(p-1)}},
 \|D^{1-s_p} I u\|_{L_{[0,t]}^{p-1} L_x^{3(p-1)}} \right\}.
\]
This function is continuous and $\tilde{Z}(0,u) = 0$. By the conclusion (\ref{st2}) of Step 2, we have
\[
 \tilde{Z}(t,u) \lesssim_{s,p} 1 + t^{\frac{5-p}{2}} \tilde{Z}^{p-1}(t,u)
 + \frac{1}{N^{\frac{5-p}{2}}} \tilde{Z}^{p-1}(t,u).
\]
By a continuity argument it is clear that there exist $N_0(s,p)$ and $\tau_0(s,p)$, such that
if $t < \tau_0$ and $N > N_0$, then $\tilde{Z}(t,u) \lesssim 1$. Plugging it back to (\ref{st2}), we finish the proof of this lemma.
\section{Proof of Almost Conservation Law of Energy}
In this section we will prove the almost conservation law of energy.
\paragraph{The Variation of the Energy}
The following computation shows the difference of the energy from time $t_1$ to time $t_2$.
\begin{eqnarray*}
 \lefteqn{E(t_2) - E(t_1) = \int_{t_1}^{t_2} \partial_t \left[\int_{\Rm^3} (\frac{|\nabla Iu|^2}{2} + \frac{|\partial_t Iu|^2}{2} + \frac{|Iu|^{p+1}}{p+1}) dx\right] dt}\\
 &=& \int_{t_1}^{t_2} \int_{\Rm^3} (\nabla Iu \cdot \partial_t \nabla Iu + \partial_{t}^2 Iu \cdot \partial_t Iu - F(Iu) \partial_t Iu) dx dt\\
 &=& \int_{t_1}^{t_2} \int_{\Rm^3} (\nabla Iu \cdot \partial_t \nabla Iu + I \Delta u \cdot \partial_t Iu  + I F(u) \cdot \partial_t Iu - F(Iu) \partial_t Iu) dx dt\\
 &=& \int_{t_1}^{t_2} \int_{\Rm^3} [\left(I F(u)-  F(Iu)\right) \cdot \partial_t Iu] dx dt.
\end{eqnarray*}
Here we use the equation (\ref{equ}).
\paragraph{The Establishment of Almost Conservation of Energy}
From the computation above we can estimate the difference by the Holder's Inequality
\[
 |E(I u(t_2)) - E(I u(t_1))| \leq \|\partial_t Iu\|_{L_J^\infty L_x^2} \|F(Iu) - I F(u)\|_{L_J^1 L_x^2}.
\]
Thus
\[
 |E(I u(t_2)) - E(I u(t_1))| \lesssim \sup_{t \in J} E (t)^{1/2} \|F(Iu) - I F(u)\|_{L_J^1 L_x^2}.
\]
The rest of the section consists of the proof of the following estimates, which immediately imply the almost conservation law.
\begin{equation} \label{esfi1}
 \|F(Iu) - F(u)\|_{L_J^1 L_x^2} \lesssim \frac{Z^p (J,u)}{N^{(5-p)/2}}.
\end{equation}
\begin{equation} \label{esfi2}
 \|F(u) - I F(u)\|_{L_J^1 L_x^2} \lesssim \frac{Z^p (J,u)}{N^{(5-p)/2}}.
\end{equation}
\paragraph{Proof of (\ref{esfi1})}
We have
\begin{eqnarray*}
\lefteqn{\|F(Iu) - F(u)\|_{L_J^1 L_x^2}
 \lesssim \|(|Iu| + |u|)^{p-1} |Iu -u|\|_{L_J^1 L_x^2}}\\
 &\lesssim& \|P_{\ll N} u\|_{L_J^{\frac{4(p-1)}{7-p}} L_x^{\frac{4(p-1)}{p-3}}}^{p-1} \cdot
\| P_{\gtrsim N} u \|_{L_J^\frac{4}{p-3} L_x^\frac{4}{5-p}}\\
&& + \| P_{\gtrsim N} u \|_{L_J^{p} L_x^{2p}}^{p-1} \cdot \|P_{\gtrsim N} u\|_{L_J^p L_x^{2p}}\\
 &\lesssim& \|D^{1-1} Iu\|_{L_J^{\frac{4(p-1)}{7-p}} L_x^{\frac{4(p-1)}{p-3}}}^{p-1} \cdot
 \frac{1}{N ^{(5-p)/2}} \| D^{1 - \frac{p-3}{2}} I u \|_{L_J^\frac{4}{p-3} L_x^\frac{4}{5-p}}\\
&& + \frac {1}{N^{(5-p)/2}} \| D^{1 - \frac{3p - 5}{2p}} Iu \|_{L_J^p L_x^{2p}} ^ p\\
 &\lesssim& \frac{Z^p (J, u)}{N^{(5-p)/2}}.
\end{eqnarray*}
Here we used the inequality
\[
 s > \frac{p-3}{2}, \frac{3p-5}{2p}.
\]
\paragraph{Proof of (\ref{esfi2})}
For the second inequality
\begin{eqnarray*}
&& \|F(u) - I F(u)\|_{L_J^1 L_x^2}\\
 &\lesssim& \|P_{\gtrsim N} F(u)\|_{L_J^1 L_x^2}\\
 &\lesssim& \|P_{\gtrsim N} F(P_{\ll N} u)\|_{L_J^1 L_x^2} +
 \|P_{\gtrsim N} [ F(u) -F(P_{\ll N} u)] \|_{L_J^1 L_x^2}\\
 &\lesssim& \|P_{\gtrsim N} F(P_{\ll N} u)\|_{L_J^1 L_x^2} +
 \| F(u) -F(P_{\ll N} u) \|_{L_J^1 L_x^2}\\
 &\lesssim& \|P_{\gtrsim N} F(P_{\ll N} u)\|_{L_J^1 L_x^2} +
 \left\| |P_{\ll N} u|^{p-1} P_{\gtrsim N} u\right\|_{L_J^1 L_x^2}\\
&& + \left\| |P_{\gtrsim N} u|^{p-1} P_{\gtrsim N} u\right\|_{L_J^1 L_x^2}\\
 &\lesssim& X_1 + X_2 + X_3.
\end{eqnarray*}
The last two terms $X_2$ and $X_3$ can be estimated in the same way as in the proof of (\ref{esfi1}),
thus we only need to consider the first term here.
\begin{eqnarray*}
 X_1 &\lesssim& \frac{1}{N^{(5-p)/2}} \|D^{\frac{5-p}{2}} F(P_{\ll N} u) \|_{L_J^1 L_x^2}\\
 &\lesssim& \frac{1}{N^{(5-p)/2}} \|P_{\ll N} u\|_{L_J^{\frac{4(p-1)}{7-p}} L_x^{\frac{4(p-1)}{p-3}}}^{p-1} \cdot
\| D^{\frac{5-p}{2}} P_{\ll N} u \|_{L_J^\frac{4}{p-3} L_x^\frac{4}{5-p}}\\
 &\lesssim& \frac{1}{N^{(5-p)/2}} \| D^{1-1} Iu\|_{L_J^{\frac{4(p-1)}{7-p}} L_x^{\frac{4(p-1)}{p-3}}}^{p-1} \cdot
 \| D^{1 - \frac{p-3}{2}} Iu\|_{L_J^\frac{4}{p-3} L_x^\frac{4}{5-p}}\\
 &\lesssim& \frac{Z^p (J,u)}{N^{(5-p)/2}}.
\end{eqnarray*}


\end{document}